\input ./style/arxiv-vmsta.cfg
\documentclass[numbers,compress,v1.0.1]{vmsta}
\usepackage{vtexbibtags}

\volume{6}
\issue{1}
\pubyear{2019}
\firstpage{41}
\lastpage{55}
\aid{VMSTA125}
\doi{10.15559/18-VMSTA125}
\articletype{research-article}

\startlocaldefs

\allowdisplaybreaks

\newtheorem{Theorem}{Theorem}[section]

\theoremstyle{definition} 
\newtheorem{Definition}{Definition}[section]
\newtheorem{ex}{Example}[section]

\endlocaldefs

\begin{document}

\begin{frontmatter}
\pretitle{Research Article}

\title{Studies on generalized Yule models}

\author{\inits{F.}\fnms{Federico}~\snm{Polito}\ead[label=e1]{federico.polito@unito.it}\orcid{0000-0003-1971-214X}}
\address{Dipartimento di Matematica, \institution{Universit\`a di Torino},
Via Carlo Alberto 10, 10123, Torino, \cny{Italy}}

\markboth{F. Polito}{Studies on generalized Yule models}

\begin{abstract}
We present a generalization of the Yule model for macroevolution in which,
for the appearance of genera, we consider point processes with the order statistics 
property, while for the growth of species we use nonlinear time-fractional
pure birth processes or a critical birth-death process. Further, in
specific cases we derive the explicit form of the distribution of the
number of species of a genus chosen uniformly at random for each time. 
Besides, we introduce a time-changed mixed Poisson process with the same
marginal distribution as that of the time-fractional Poisson process.
\end{abstract}

\received{\sday{16} \smonth{7} \syear{2018}}
\revised{\sday{10} \smonth{11} \syear{2018}}
\accepted{\sday{16} \smonth{11} \syear{2018}}
\publishedonline{\sday{3} \smonth{12} \syear{2018}}
\end{frontmatter}

\section{Introduction}

In 1925, Udny Yule published the paper \cite{yule1925ii} in which he
described a possible model for macroevolution.\index{macroevolution} Genera (species
grouped by similar characteristics) appear in the system at random
times\index{random ! times} according to a linear birth process.
Each genus is initially composed
by a single species.\index{single species} As soon as the genus appears, an independent
linear birth process modelling the evolution of the species belonging
to it, starts. The initial species thus generates offsprings
(representing related species) with constant individual intensities.
The model is now classical
and it is usually named after him.
The Yule model\index{Yule model} therefore admits the existance of two independent and
superimposed mechanisms of evolution, one for
genera and one for species, both at possibly different constant intensities.

One of the most interesting characteristics the model exhibits is its
intrinsic ``preferential attachment\index{preferential attachment}'' mechanism.
And this is exactly implied by the presence of the two distinct and
superimposed random growths
for the appearance of genera and of species.

The preferential attachment\index{preferential attachment} mechanism is a fundamental ingredient of
many modern models of random graphs\index{random ! graphs} growth.
The literature is vast and covers many fields. We recall here only
some relevant papers and books
\cite{Athreya,MR2091634,berger2014,DiameterBAmodel,Borgs:2007:FME:1250790.1250812,Remco92,CPC:9224998,Faloutsos:1999,GaravagliaHofstad,Kleinberg:1999,redner,PhysRevLett.85.4629,MR3163720,yang,hofstad2016} leaving the reader
the possibility of widening the number of sources by looking at the
references cited therein.

The Barab\'asi--Albert model of random graph\index{random ! graphs} growth is by far the most
famous example
of a stochastic process based on preferential attachment\index{preferential attachment} \cite
{MR2091634,MR1824277}.
In \cite{PachonPolitoSacerdote} we discuss the relationships between
the Barab\'asi--Albert graph, the Yule model\index{Yule model} and a third model
based on preferential attachment\index{preferential attachment} introduced by Herbert Simon in 1955
\cite{lansky1980generalization,SIMON1955,MR0130733}.
In \cite{pachon2016continuous} we have further analyzed and described
the exact relation between the Barab\'asi--Albert model and
the Yule model.\index{Yule model} Briefly, the finite-dimensional distributions of the
degree of a vertex in the Barab\'asi--Albert model
converges to the finite-dimensional distributions of the number of
individuals in a Yule process\index{Yule process} with initial population size equal
to the number $m$ of attached edges in each time step. This further
entails that the asymptotic degree distribution\index{degree distribution} of
a vertex chosen uniformly at random in the Barab\'asi--Albert model
coincides with the asymptotic distribution
of the number of species belonging to a genus chosen uniformly at
random from the Yule model\index{Yule model} with an initial number $m$ of species.
This result suggests that asymptotic models similar to the Yule model\index{Yule model}
can be linked to different preferential attachment random
graph\index{preferential attachment random graph} processes in discrete time.

Therefore, a direct analysis of the Yule model,\index{Yule model} a model in continuous
time and hence possessing a greater mathematical
treatability, has the potential to uncover important aspects and
characteristics of the discrete-time model to which it is related.

In \cite{lansky2014role} and \cite{MR3562427} we have started a study
of macroevolutionary models similar to the classical Yule model where
the process governing the appearance of genera is left unchanged,
while those describing the growth of species
account for more realistic features. Specifically, in \cite
{lansky2014role} we have generalized the latter allowing the
possibility of extinction of species while in \cite{MR3562427} we have
studied the effect of a slowly-decaying memory
by considering a fractional nonlinear birth process. Notice that, by
suitably specializing the nonlinear rates, other peculiar
behaviours such as saturation or logistic growth may be observed.

In this paper we proceed with the analysis by looking at modifications
of the classical Yule model in which the appearance of genera follows
a different dynamics. We will show that a change in the dynamics of
genera will lead to radical changes in the model.
This is a preliminary step before aiming at deriving models of random
graphs\index{random ! graphs} with different features than those graphs
connected to the Yule model.\index{Yule model}
In the following we consider the class of mixed Poisson processes
time-changed by means of a deterministic
function. The rationale which justifies this choice will be clear in
Section \ref{main} where we state the main results.
Section \ref{prel} contains the mathematical background necessary to
develop the results presented later.
In particular we will connect a member of the class of the suitably
time-changed mixed Poisson processes with the
time-fractional Poisson process, a non-Markov renewal process governed
by a time-fractional difference-differential equations
involving the Caputo--D\v{z}rba\v{s}jan derivative.

\section{Preliminaries}\label{prel}

We consider two different classes of point processes, namely the
time-fractional Poisson processes (shortly tfPp\index{tfPp}) and an extension of
the mixed Poisson processes, i.e.\ the mixed Poisson processes up to a
(deterministic) time transformation (mPp-utt in the following).
We will analyze briefly their properties and state a result linking
the two classes.
Besides, we will introduce the definitions and mathematical tools
needed to understand Section \ref{main}.

\subsection{Time-fractional Poisson processes (tfPp\index{tfPp})}

The tfPp\index{tfPp} has been introduced in the literature in \cite
{MR1604710,MR1910034} (see also Laskin's paper \cite{MR2007003}).
We show here the construction by means of random time-change with an
inverse stable subordinator\index{stable subordinator} \cite{MR2835248}.
Alternatively, the tfPp\index{tfPp} can be defined as a specific renewal process
(see e.g.\ \cite{MR2535014,MR2650778,MR2120631},
and see \cite{MR2835248} for the proof of the equivalence between the
two constructions).

Let us consider a homogeneous Poisson process\index{homogeneous Poisson process} $(\widetilde{N}_t)_{t
\ge0}$ of parameter $\lambda>0$ and an independent inverse stable
subordinator,\index{stable subordinator} that is a one-dimensional time-continuous stochastic and
non-Markov process defined as follows.
Consider the subordinator $(D_t)_{t \ge0}$ with L\'evy measure $\nu
(\textup{d}x)=[\alpha/\varGamma(1-\alpha)]x^{-1-\alpha}$,
$\alpha\in(0,1)$, and define its stochastic inverse $(E_t)_{t \ge
0}$ as the first time at which it exceeds a given threshold, i.e.
%
\begin{align}
E_t=\inf\{ s \colon D_s > t \}, \quad t \ge0.
\end{align}
Now, consider the time-changed point process $\mathcal{N}=(\mathcal
{N}_t)_{t \ge0}=(\widetilde{N}_{E_t})_{t \ge0}$.
The process $\mathcal{N}$ is called tfPp of parameters $\lambda$ and
$\alpha$.

Many properties are known for the tfPp.\index{tfPp} Let us review some of them.
The mar\-ginal~probability distribution of the process
generalizes the Poisson distribution (for $\alpha\to1$) and can be
written as \cite{MR2535014,MR2650778}
%
\begin{align}
\label{colla} \mathbb{P}(\mathcal{N}_t=k) = \bigl(\lambda
t^\alpha\bigr)^k E_{\alpha,\alpha
k+1}^{k+1} \bigl(-
\lambda t^\alpha\bigr), \quad k \ge0,\ t \ge0,
\end{align}
where
%
\begin{align}
E^\gamma_{\nu, \beta}(z) = \sum_{r=0}^\infty
\frac{(\gamma)_r z^r}{r!
\varGamma(\nu r+\beta)}, \qquad(\gamma)_r = \frac{\varGamma(\gamma+r)}{\varGamma(\gamma)}, \quad z \in
\mathbb{C},
\end{align}
is the Prabhakar function\index{Prabhakar function} \cite{MR0293349} for complex parameters $\nu
,\beta,\gamma$, with $\Re(\nu)>0$.
It is interesting to note that the mean value of the process grows in
time less than linearly for each allowed value of $\alpha$,
%
\begin{align}
\mathbb{E}(\mathcal{N}_t) = \lambda t^\alpha/ \varGamma(
\alpha+1), \quad t \ge0,
\end{align}
and that the variance can be written as
%
\begin{align}
\mathbb{V}\text{ar}(\mathcal{N}_t) = \mathbb{E}(
\mathcal{N}_t) + \frac
{(\lambda t^\alpha)^2}{\alpha} \biggl( \frac{1}{\varGamma(2\alpha)} -
\frac{1}{\alpha\varGamma(\alpha)^2} \biggr),
\end{align}
thus highlighting the overdispersion of the process. Furthermore, the
probability generating function\index{probability generating function} reads
%
\begin{align}
G(u,t) = E_\alpha\bigl(\lambda(u-1)t^\alpha\bigr), \quad|u|
\le1,
\end{align}
where $E_\alpha(z) = E_{\alpha,1}^1(z)$ is the classical
Mittag-Leffler function.

Considering the renewal nature of $\mathcal{N}$, and calling $U_j$,
$j\ge1$, the random inter-arrival times between
the $(j-1)$-th and the $j$-th event, it is possible to give an
explicit expression for the common probability density function.
Indeed,
%
\begin{align}
f_{U_j}(t) = \lambda t^{\alpha-1} E_{\alpha,\alpha}\bigl(-\lambda
t^\alpha\bigr) \mathbb{I}_{\mathbb{R_+}}(t),
\end{align}
where $E_{\alpha,\alpha}(z)=E^1_{\alpha,\alpha}(z)$ is the generalized
Mittag-Leffler function. By analyzing the density,
its slowly decaying right tail (it is actually an ultimately monotone
regularly varying function of order $-\alpha-1$)
and the asymptote in zero, the clusterization of events in time
appears evident.

Lastly, let us recall the direct relation linking the tfPp\index{tfPp} with
fractional calculus: the probability distribution of the tfPp
solves a specific difference-differential equation in which the time
derivative appearing in the difference-differential equations
related to the homogeneous Poisson process\index{homogeneous Poisson process} is replaced by a fractional
derivative of the Caputo--D\v{z}rba\v{s}jan type.
Regarding this, for $ n\in\mathbb{N}$, denote by $\textup{AC}^{m}
[a,b]$ the
space of real-valued functions with
continuous derivatives up to order $m-1$
such that the $(m-1)$-th derivative belongs to the space of absolutely
continuous functions
$\textup{AC}[a,b]$. In other words,
%
\begin{equation}
\textup{AC}^m [a,b ] = \biggl\{ f\colon [a,b ] \mapsto \mathbb{R}
\colon\frac{\textup{d}^{m-1}}{\textup{d}x^{m-1}} f ( x ) \in\textup{AC} [a,b] \biggr\} .
\end{equation}
Then, for $\alpha>0$, $m = \lceil\alpha\rceil$, and $f \in
AC^m[a,b]$, the Caputo--D\v{z}rba\v{s}jan derivative of order
$\alpha>0$ is defined as
%
\begin{equation}
\label{Capu} {}^C D^{\alpha}_{a^+}f(t)=
\frac{1%
}{\varGamma(m-\alpha)}\int_a^{t}(t-s)^{m-1-\alpha}
\frac{\textup
{d}^m}{\textup{d}s^m}f(s) \textup{d}s.
\end{equation}
Let us now denote the state probabilities $\mathbb{P}(\mathcal{N}_t=k)$
of the tfPp\index{tfPp} by $p_k(t)$, $k \ge0$, $t \ge0$.
Then the probabilities $p_k(t)$ satisfy the equations
%
\begin{align}
\label{eqeq} {}^C D^{\alpha}_{0^+}
p_k(t) = -\lambda p_k(t) + \lambda p_{k-1}(t),
\quad k \ge0,
\end{align}
where we consider $p_{-1}(t)$ being equal to zero.
In Example \ref{equiv}, the tfPp\index{tfPp} will be compared with a member of the
class of the mPp-utt.

\subsection{Mixed Poisson processes up to a time transformation
(mPp-utt)}\label{mpp-utt}

We start describing mixed Poisson processes (mPp), first introduced by
J.\ Dubourdieu in 1938 \cite{dubourdieu1938remarques}.
For full details the reader can refer to the monograph by J.\ Grandell
\cite{MR1463943}.

Consider a unit-rate homogeneous Poisson process $N=(N_t)_{t \ge0}$.
A point process $(\widetilde{M}_t)_{t \ge0}$ is an mPp if and only if
$\widetilde{M}_t = N_{Wt}$ in distribution, where
$W$ is an almost surely non-negative random variable\index{random ! variable}
independent of $N$.
Common choices for the mixing random variable $W$\index{random ! variable} are the Gamma
distribution, leading to the P\'olya process, or the uniform
distribution on $[0,c)$, $c \in\mathbb{R}_+$.

Clearly, if $W$ is degenerate on $w$, then an mPp coincides with a
homogeneous Poisson process\index{homogeneous Poisson process} of rate $w$.

An mPp is characterized by the so-called property $P$ which means that,
conditional on $\widetilde{M}_t-\widetilde{M}_0=k$,
the random jump times $\{ t_1,t_2 \dots, t_k \}$ are distributed as
the order statistics
of $k$ iid uniform random variables\index{random ! variable} on $[0,t]$ \cite{MR531764}.

This result, first appeared in \cite{MR0154323}, can be further
extended considering a deterministic time-change, leading
to the class of point processes with the \textit{OS} property (order
statistics property)
\cite{MR0438474,MR531764}: a point process $(K_t)_{t \ge0}$ with unit
steps is said to have the
\textit{OS} property if and only if, conditional on $K_t-K_0=k$, the random
jump times $\{ t_1,t_2 \dots, t_k \}$ are distributed as the order statistics
of $k$ iid random variables supported on $[0,t]$ with distribution
function $F_t(x) = q(x)/q(t)$, where $q(t) = \mathbb{E}(K_t)-\mathbb
{E}(K_0)$ is continuous and non-decreasing. In this respect, property
$P$ is also called uniform \textit{OS} property.

Notably, K.S.\ Crump proved that point processes with the \textit{OS}
property\index{OS property} are Markovian (see \cite{MR0438474}, Theorem 2).

Taking into account the results presented in \cite{MR0438474,MR531764}, and \cite{MR644418},
we recall the following theorem due to P.D.\ Feigin \cite{MR531764}:
%
\begin{Theorem}
\label{f}
Let $M$ be a point process with the $OS$ property relative to the
distribution function
$F_t(x) = q(x)/q(t)$, where $q(t) = \mathbb{E}(M_t) -\mathbb{E}(M_0)$
is a
continuous and non-decreasing function. Then there exists a unit-rate
homogeneous Poisson
process $N$ and an independent non-negative random variable\index{random ! variable} $W$
defined on the same probability
space, such that $M_t = N_{Wq(t)}$ almost surely.
\end{Theorem}

Notice also that Theorem \ref{f} implies $\mathbb{E}(W)=1$.
Furthermore, the above theorem does not exclude the case of bounded
$q$, that is when $\lim_{t \rightarrow\infty}
q(t) = \gamma< \infty$. Processes in that subclass are usually called
mixed sample processes.
To gain more insight on them, the reader can consult \cite{MR644418},
Section 2, in which an interesting
example is described (see also \cite{MR789355}).

It seems clear that the class of point processes with the \textit{OS} property\index{OS property} contains
that of Mixed Poisson processes up to the time transformation $q(\cdot
)$ (mPp-utt).
We will consider in the following only the subclass of mPp-utt.

Let us first present a didactic example of a member of the class of
mPp-utt, the Yule process.\index{Yule process} Being an mPp-utt, the Yule process
exhibits\index{Yule process} the \textit{OS} property.\index{OS property} The reader may refer to \cite{MR0438474}
for more details.

\begin{ex}
Let $M$ be a Yule process starting\index{Yule process} with a single individual,
shifted downwards by one and with individual splitting rate $\lambda$.
Let $N$ be a unit-rate homogeneous Poisson process and let $W$ be
independent of $N$ and exponentially distributed with mean one.
Set $q(t) = e^{\lambda t}-1$. Then, we construct the mPp-utt
representation of $M$, i.e.\
$M_t = N_{Wq(t)}$, $t \ge0$. Note that the state space of $N_{Wq(t)}$
is $\{0,1,\dots\}$. The distribution of $M$ can be derived easily
by conditioning:
%
\begin{align}
\mathbb{P} (M_t = k) & = \int_0^\infty
\frac{[w(e^{\lambda
t}-1)]^k}{k!} e^{-w(e^{\lambda t}-1)} e^{-w} \textup{d}w\notag
\\
& = \frac{(e^{\lambda t}-1)^k}{k!} \int_0^\infty
w^k e^{-w e^{\lambda
t}} \textup{d}w \notag
\\
& = e^{-\lambda t}\bigl(1-e^{-\lambda t}\bigr)^k, \quad k \ge0.
\end{align}
The Yule process\index{Yule process} has the \textit{OS} property with $F_t = \frac{e^{\lambda
x}-1}{e^{\lambda t}-1}$, $x \in[0,t]$.
The \textit{OS} property of the Yule process\index{Yule process} has been implicitly used in many
papers, starting from the seminal paper by Yule
\cite{yule1925ii} (see also \cite
{chan2003stochastically,lansky2014role,MR3562427})
\end{ex}

In the following example we define a specific mPp-utt which is
connected with the tfPp\index{tfPp} by means of its marginal
distribution.\index{marginal distribution}

\begin{ex}
\label{equiv}
Let $M = (M_t)_{t \ge0}$ be such that $M_t = N_{Wq(t)}$ where $q(t) =
\lambda t^\nu/\varGamma(1+\nu)$,
$\lambda> 0$, $\nu\in(0,1)$, and $W$ is a unit-mean non-negative
random variable\index{random ! variable}
with probability density function
%
\begin{align}
f_{W}(w) =  \phi \biggl(-\nu,1-\nu; \frac{-w}{\varGamma(1+\nu)} \biggr)
\bigg/ \varGamma(1+\nu),  \quad w \in\mathbb{R}_+.
\end{align}
The above density is written in terms of the Wright function (see \cite
{MR2218073} for details)
%
\begin{align}
\phi(\alpha,\beta;z) = \sum_{r=0}^\infty
\frac{z^r}{r! \varGamma(\alpha
r+\beta)}, \quad\alpha,\beta,z \in\mathbb{C}, \ \Re(\alpha) > -1.
\end{align}

Plainly, $\mathbb{E}(M_t) = \lambda t^\nu/\varGamma(1+\nu)$. Now, let us
derive the marginal distribution\index{marginal distribution} of the mPp-utt $M$:
%
\begin{align}
\mathbb{P} (M_t = k) & = \int_0^\infty
\mathbb{P} (N_{wq(t)}=k) f_{W}(w) \textup{d} w\notag
\\
& = \int_0^\infty\frac{ ( w\frac{\lambda t^\nu}{\varGamma(1+\nu)}
 )^k}{k!}
e^{-w\frac{\lambda t^\nu}{\varGamma(1+\nu)}} \phi \biggl(-\nu,1-\nu; \frac{-w}{\varGamma(1+\nu)} \biggr)
\frac{\textup
{d}w}{\varGamma(1+\nu)}.
\end{align}
By letting $\xi= wt^\nu/\varGamma(1+\nu)$ we have
%
\begin{align}
\mathbb{P} (M_t = k) & =\int_0^\infty
\frac{(\lambda\xi)^k}{k!} e^{-\lambda\xi} t^{-\nu}\phi\bigl(-\nu,1-\nu;-\xi
t^{-\nu}\bigr) \textup{d} \xi\notag
\\
& = \mathbb{P} (\widetilde{N}_{E_t} = k) = \mathbb{P}(\mathcal
{N}_t=k).
\end{align}
This last step is justified by the time-change construction of the tfPp\index{tfPp}
and by the fact that the marginal density function of the
inverse stable subordinator\index{stable subordinator} $(E_t)_{t \ge0}$ is exactly $f_{E_t}(\xi
)=t^{-\nu}\phi(-\nu,1-\nu;-\xi t^{-\nu})$, $\xi\in\mathbb{R}_+$ (see
e.g.\ \cite{MR2845910}, Section 2).

It is worthy of note that the tfPp\index{tfPp} $\mathcal{N}$ and the mPp-utt $M$
share the same marginal distribution.\index{marginal distribution} This entails
that the probabilities $\mathbb{P} (M_t=k)$, $k \ge0$, solve the
equations~\eqref{eqeq}.
\end{ex}

\section{Generalized Yule model}\label{main}

We proceed now to the analysis of a generalization of the Yule model\index{Yule model}
in the sense we have anticipated in the introductory section.
The focus here is to construct a model in which the arrival in time of
genera is driven by an mPp-utt and the
process describing the evolution of species for each different genus
is a tfPp\index{tfPp} or a nonlinear time-fractional
pure birth process (see \cite{MR3562427,MR2730651}).
Hence what we drop here is the deterministic constant intensity
assumption for genera evolution.
The genera arrival is instead described by random intensities.\index{random ! intensities}
For the sake of clarity, before describing the generalization of the
Yule model,\index{Yule model} let us recall the definition of the
nonlinear time-fractional pure birth process. Analogously as for the
tfPp,\index{tfPp} the construction
of the nonlinear time-fractional pure birth process is by time-change
with the inverse stable subordinator\index{stable subordinator} $(E_t)_{t \ge0}$.
Thus, consider the nonlinear pure birth process $(Y_t)_{t \ge0}$,
starting with a single progenitor, with nonlinear rates
$\lambda_k > 0$, $k \ge1$, and being independent of $(E_t)_{t \ge
0}$. The time-changed process
$\mathcal{Y}=(\mathcal{Y}_t)_{t \ge0} = (Y_{E_t})_{t \ge0}$ is called
a nonlinear time-fractional pure birth process. Interestingly enough,
when $\lambda_k = \lambda$ for each $k \ge1$, the
nonlinear time-fractional pure birth process reduces to the tfPp\index{tfPp} of
parameters $\lambda$ and $\nu$, shifted upwards by one.

Let us now consider the following model.

\begin{Definition}[Generalized Yule model]
The generalized Yule model represents the growth of a population which
evolves according to:
\begin{enumerate}
\item Genera (each initially with a single species\index{single species}) appear following an
mPp-utt\break $(M_t)_{t\geq0}$.
\item When a new genus appears a copy of $\mathcal{Y}$ starts. The
copies are independent one of another and of the
mPp-utt. Each copy models the evolution of species belonging to the
same genus.
\end{enumerate}
\end{Definition}

Then, for each time $t \in\mathbb{R}_+$ we define the random variable\index{random ! variable}
$\mbox{}_t \mathfrak{N}$ measuring the number of species
belonging to a genus chosen uniformly at random. With respect to the
classical Yule model this random variable\index{random ! variable} is linked to
the degree distribution\index{degree distribution} of a vertex chosen uniformly at random in the
Barab\'asi--Albert model \cite{pachon2016continuous}.
Our aim is to investigate the distribution of $\mbox{}_t \mathfrak{N}$ for
the generalized Yule model.
To do so, it is enough to condition on the random creation time $T$ of
the selected genus, obtaining
%
\begin{align}
\mathbb{P}(\mbox{}_t\mathfrak{N} = k) = \mathbb{E}_T
\mathbb{P} (\mathcal {Y}_t=k |\mathcal{Y}_T=1), \quad k
\ge1.
\end{align}
Notice that, due to the $OS$ property satisfied by the considered
mPp-utt, the distribution function of $T$ is $F_t(\cdot)$
(see Section \ref{mpp-utt}).

We specialize now the model by choosing the process of Example \ref
{equiv} for the random arrival of genera.

\begin{figure}[t!]
\includegraphics{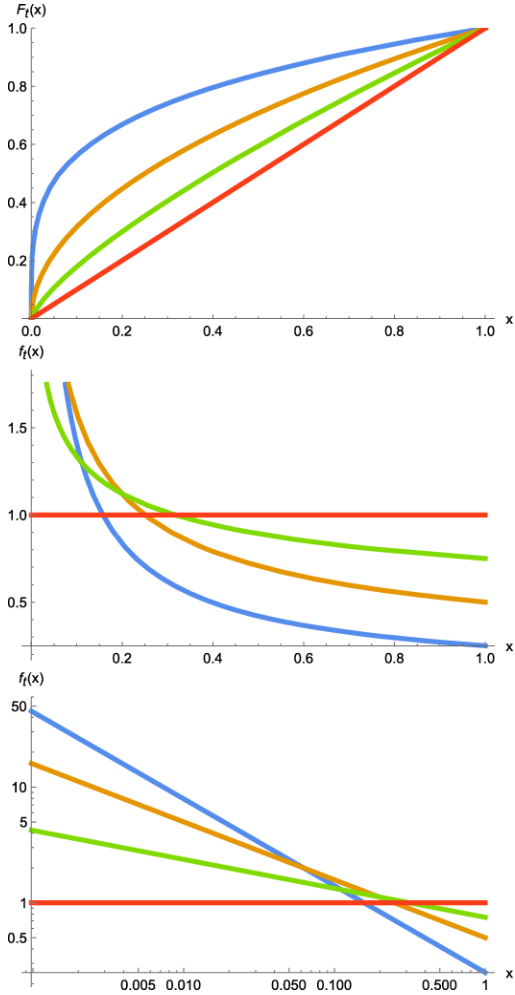}
\caption{Distribution function \eqref{ad} (top) and
density function \eqref{bd}
in linear plot (middle) and loglog plot (bottom). The parameter $\nu
$ is set to
$\nu=(1/4,1/2,3/4,1) = (\text{blue},\text{orange},\text{green},\text{red})$
and $t=1$. Note how, for $\nu\in(0,1)$ (in contrast with the
classical case $\nu=1$),
the concentration of probability mass near zero makes the appearance
of genera more likely to occur in the very early
evolution of the process}\label{distri}
\end{figure}

In this case the
distribution function of $T$ reads
%
\begin{align}
\label{ad} F_t(x) = \biggl( \frac{x}{t}
\biggr)^\nu, \quad x \in[0,t],
\end{align}
with density
%
\begin{align}
\label{bd} f_t(x) = \frac{\nu x^{\nu-1}}{t^\nu}, \quad x \in[0,t].
\end{align}
Figure \ref{distri} shows the shapes of the distribution function \eqref
{ad} and the density function \eqref{bd}
for different values of the characterizing parameter $\nu$. Notice the
rather different behaviour for values
of $\nu$ strictly less than 1.

Regarding the evolution of the number of species for each genus, the
fractional exponent of the process $\mathcal{Y}$
will be denoted by $\beta\in(0,1)$.
We suppose now that the nonlinear rates of $\mathcal{Y}$ are all
different and recall that in this case
%
\begin{align}
\label{nonli} \mathbb{P}(\mathcal{Y}_t = k) = \prod
_{j=1}^{k-1} \lambda_j \sum
_{m=1}^k \frac{E_\beta(-\lambda_m t^\beta)}{\prod_{l=1,l\ne m}^k (\lambda
_l-\lambda_m)}, \quad k \ge1,
\end{align}
with the convention that empty products equal unity. We obtain
%
\begin{align}
\mathbb{P}(\mbox{}_t\mathfrak{N} = k) & = \int_0^t
\mathbb{P}(\mathcal {Y}_{t-x}=k) \frac{\nu x^{\nu-1}}{t^\nu} \textup{d} x\notag
\\
& = \frac{\nu}{t^\nu}\prod_{j=1}^{k-1}
\lambda_j \sum_{m=1}^k
\frac
{1}{\prod_{l=1,l\ne m}^k (\lambda_l-\lambda_m)} \int_0^t E_\beta
\bigl[-\lambda_m(t-x)\bigr] x^{\nu-1} \textup{d}x.
\end{align}
Now we make use of Corollary 2.3 of \cite{MR2033353} and arrive at
%
\begin{align}
\label{formulona} \mathbb{P}(\mbox{}_t\mathfrak{N} = k) = \varGamma(\nu+1)
\prod_{j=1}^{k-1} \lambda_j \sum
_{m=1}^k \frac{E_{\beta,\nu+1}(-\lambda_m t^\beta)}{\prod_{l=1,l\ne m}^k
(\lambda_l-\lambda_m)}, \quad k \ge1.
\end{align}

\begin{figure}[t!]
\includegraphics{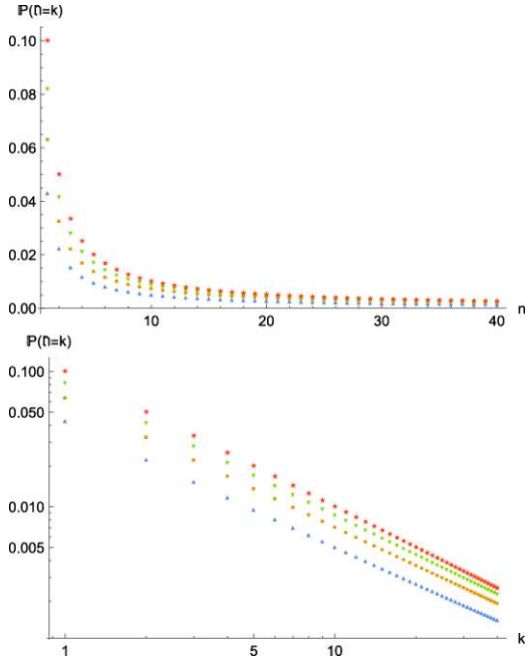}
\caption{The probability mass function of $\mbox{}_t\mathfrak
{N}$ in the linear rate case
(formula \eqref{formulina}). Top: linear plot, bottom: loglog plot,
with $t=10$, $\beta=1$ (which
corresponds to classical linear birth processes for species
evolution), $\lambda=1$, and
$\nu=(0.4, 0.6, 0.8, 1)$ from bottom to top}\label{pmf}
\end{figure}

A special case of interest is when the rates are linear, $\lambda
_k=\lambda k$, $k \ge1$. In this case, from
\eqref{formulona} we obtain easily the following probabilities:
%
\begin{align}
\label{formulina} \mathbb{P}(\mbox{}_t\mathfrak{N} = k) = \varGamma(\nu+1)
\sum_{j=1}^k \binom
{k-1}{j-1}(-1)^{j-1}
E_{\beta, \nu+1}\bigl(-\lambda j t^\beta\bigr), \quad k \ge1.
\end{align}
Figure \ref{pmf} shows how the above probability mass function changes
with respect to parameter $\nu$,
taking a constant $\beta=1$, that is, considering a classical behaviour
for species.

Recalling that in the linear rates case $\mathbb{E}\mathcal
{Y}_t=E_\beta(\lambda t^\beta)$
and $\mathbb{E}\mathcal{Y}_t^2= 2E_\beta(2\lambda t^\beta) - E_\beta
(\lambda t^\beta)$, we derive
the first two moments for the random variable\index{random ! variable} $\mbox{}_t\mathfrak{N}$ and
its variance:
%
\begin{align}
\mathbb{E}\;\mbox{}_t\mathfrak{N}&=\frac{\nu}{t^\nu} \int
_0^t E_\beta\bigl[\lambda
(t-x)^\beta\bigr] x^{\nu-1}\mathrm{d}x = \varGamma(\nu+1)
E_{\beta,\nu+1}\bigl(\lambda t^\beta\bigr),
\\
\mathbb{E}\;\mbox{}_t\mathfrak{N}^2 &= \frac{\nu}{t^\nu}
\int_0^t\mathbb {E}\mathcal{Y}_{t-x}^2
x^{\nu-1} \mathrm{d}x\notag
\\
& = 2\varGamma(\nu+1) E_{\beta,\nu+1}\bigl(2\lambda t^\beta\bigr) -
\varGamma(\nu+1) E_{\beta,\nu+1}\bigl(\lambda t^\beta\bigr),
\\
\mathbb{V}\text{ar}\; \mbox{}_t\mathfrak{N} &= 2\varGamma(\nu+1)
E_{\beta
,\nu+1}\bigl(2\lambda t^\beta\bigr)\notag
\\
&\quad  - \varGamma(\nu+1)E_{\beta,\nu+1}\bigl(\lambda t^\beta\bigr)
\bigl(1+ \varGamma(\nu +1)E_{\beta,\nu+1}\bigl(\lambda t^\beta\bigr)
\bigr).
\end{align}

When the nonlinear rates are actually constant and all equal (i.e.\
$\lambda_k=\lambda$, $\forall k \ge1$) we cannot make use of a
specialized form of formula
\eqref{nonli}. In this case however, the nonlinear time-fractional
pure birth process reduces to the tfPp suitably shifted upwards by one.
Hence, recalling formula \eqref{colla}, the distribution of
$\mbox{}_t\mathfrak{N}$
can be written by conditioning as
%
\begin{align}
\mathbb{P}(\mbox{}_t\mathfrak{N} = k) & = \mathbb{E}_T
\mathbb{P} (\mathcal {N}_t+1=k |\mathcal{N}_T+1=1\xch{)}{),}\notag
\\
& = \int_0^t \mathbb{P}(\mathcal{N}_{t-x}=k-1)
\;\frac{\nu x^{\nu
-1}}{t^\nu} \textup{d} x, \quad k \ge1.
\end{align}
Then we have
%
\begin{align}
\mathbb{P}(\mbox{}_t\mathfrak{N} = k) & = \int_0^t
\bigl[\lambda(t-x)^\beta \bigr]^{k-1} E_{\beta,\beta(k-1)+1}^k
\bigl[-\lambda(t-x)^\beta\bigr] \, \frac{\nu x^{\nu-1}}{t^\nu} \textup{d} x\notag
\\
& = \frac{\nu\lambda^{k-1}}{t^\nu} \int_0^t
(t-x)^{\beta k -\beta} E_{\beta,\beta k -\beta+1}^k \bigl[-\lambda(t-x)^\beta
\bigr] x^{\nu-1} \textup{d}x.
\end{align}
The above integral is known and can be calculated by using Corollary
2.3 of \cite{MR2033353}. We finally obtain
%
\begin{align}
\label{colla2} \mathbb{P}(\mbox{}_t\mathfrak{N} = k) = \varGamma(\nu+1)
\lambda^{k-1} t^{\beta k -\beta} E_{\beta,\beta k - \beta+ \nu+ 1}^k \bigl(-
\lambda t^\beta\bigr), \quad k \ge1.
\end{align}

In the classical Yule model, $\lim_{t \rightarrow\infty} \mbox{}_t\mathfrak
{N} = \mathfrak{N}$, where $\mathfrak{N}$ is a
non-degenerate limiting random variable.\index{random ! variable} The distribution of $\mathfrak
{N}$ is known and is called Yule--Simon distribution.
Its main feature is the characteristic right tail which slowly decays
as a power-law.
In our cases, however, the random variable $\mbox{}_t\mathfrak{N}$\index{random ! variable} has a
different behaviour at $\infty$.
This can be observed by considering the asymptotic expansion of the
Prabhakar function\index{Prabhakar function} \cite{MR3709831}. We have
from formula \eqref{colla2}, for $t \rightarrow\infty$,
%
\begin{align}
\mathbb{P}(\mbox{}_t\mathfrak{N} = k) \sim\frac{\varGamma(\nu+1)}{\lambda}
\frac{t^{-\beta}}{\varGamma(-\beta k)} \longrightarrow0,
\end{align}
for each finite value of $k \ge1$.

\subsection{A critical macroevolutionary model with species deletion}

\begin{figure}
\includegraphics{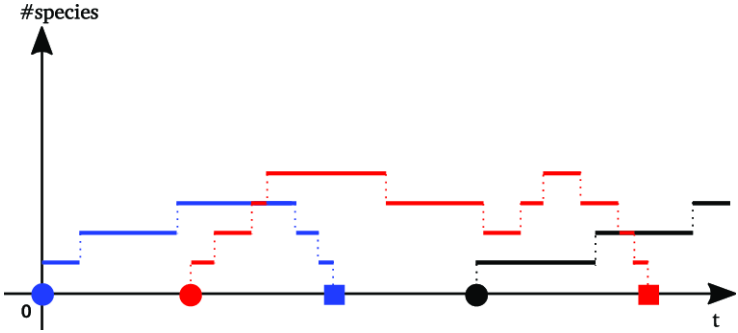}
\caption{A possible realization of the superimposed
processes counting the
number of species belonging to each existent genus (shown in different colours).
Birth of genera, governed
by the mPp-utt of Example \ref{equiv}, is represented
in the figure by circles while their extinction time is marked by squares}\label{figfig}
\end{figure}

We introduce here a model of macroevolution\index{macroevolution} in which the possibility
of extinction of genera is
taken into consideration. To achieve this, the species dynamics is
described by independent critical
birth-death processes (of parameter $\lambda>0$), each starting with a
single species,\index{single species} while the genera appearance follows the mPp-utt of
Example \ref{equiv}.
Figure \ref{figfig} shows a possible realization of the superimposed
processes counting the
number of species belonging to each existent genus. Extinction of
genera is represented by squares
while their births by circles.

Recurring to the integral representation of the Gauss
hypergeometric function,
%
\begin{align}
\mbox{}_2F_1(a,b;c;z) = \frac{\varGamma(c)}{\varGamma(b)\varGamma(c-b)} \int
_0^1 y^{b-1} (1-y)^{c-b-1}
(1-yz)^{-a} \mathrm{d} y,
\end{align}
for $c>b>0$, we derive the exact form of the transient distribution of
$\mbox{}_t \mathfrak{N}$:
%
\begin{align}
\label{hypp1} \mathbb{P}(\mbox{}_t \mathfrak{N} = 0) & =
\frac{\nu}{t^\nu} \int_0^t x^{\nu
-1}
\frac{\lambda(t-x)}{1+\lambda(t-x)} \mathrm{d}x\notag
\\
& = \frac{\lambda t}{\nu+1} \mbox{}_2F_1 (1,2;\nu+2;-\lambda t),
\end{align}
and
%
\begin{align}
\label{hypp2} \mathbb{P}(\mbox{}_t \mathfrak{N} = k) & =
\frac{\nu}{t^\nu} \int_0^t x^{\nu
-1}
\frac{\lambda^{k-1} (t-x)^{k-1}}{(1+\lambda(t-x))^{k+1}} \mathrm {d}x\notag
\\
& = (\lambda t)^{k-1} \frac{\varGamma(k)\varGamma(\nu+1)}{\varGamma(\nu+k)} \mbox{}_2F_1(k+1,k;
\nu+k;-\lambda t), \quad k \ge1.
\end{align}

%
\begin{figure}[t!]
\includegraphics{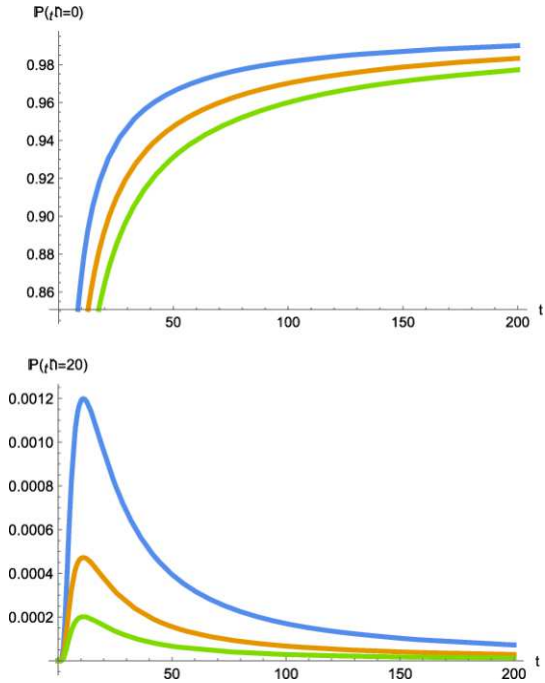}
\caption{The probabilities \eqref{hypp1} (top) and \eqref
{hypp2} (bottom, $k=20$) drawn with
respect to time, with $\lambda=1$ and $\nu=(0.2,0.5,0.8)=(\text
{blue},\text{orange},\text{green})$. Note how the probability of
selecting uniformly at random an extinct genus increases in time}\label{hyp}
\end{figure}

\begin{figure}[t!]
\includegraphics{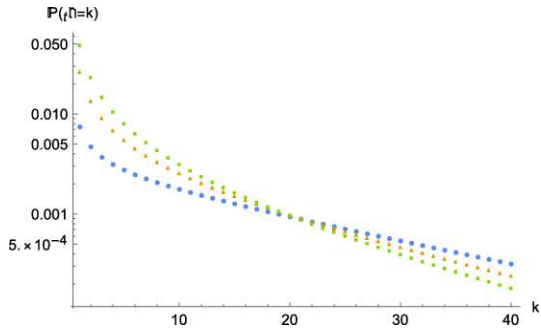}
\caption{The probability mass function \eqref{hypp2}
depicted for $\nu=(0.1,0.5,1)=(\text{blue},
\text{orange},\break \text{green})$ with $\lambda=1$}\label{puri}
\end{figure}

In Figure \ref{hyp} the above probabilities (for $k=20$) are pictured
with respect to time. The probability
mass function concentrates on zero for $t \to\infty$ as it should be.
Notably, it exhibits an exponential
tail (see in Figure \ref{puri}), differently from the case without
deletion (for example compare it with \eqref{formulina}, see Figure \ref{pmf}).
The derivation of the moments of the random variable\index{random ! variable} $\mbox{}_t\mathfrak{N}$
is simpler in this model. Recalling
that each species process has mean 1 and $\mathbb{E}\mathcal
{Y}_t^2=2\lambda t+1$ we obtain
that the expectation of $\mbox{}_t \mathfrak{N}$ is also 1 and that
%
\begin{gather}
\mathbb{E} \; \mbox{}_t \mathfrak{N}^2 = \frac{\nu}{t^\nu}
\int_0^t \bigl(2\lambda(t-x)+1\bigr)
x^{\nu-1}\mathrm{d}x = \frac{2\lambda t}{\nu+1} +1,\\
\mathbb{V}\text{ar} \; \mbox{}_t \mathfrak{N} = \frac{2\lambda t}{\nu+1}.
\end{gather}


\begin{funding}
F.\ Polito has been partially supported by the projects \emph{Memory in
Evolving Graphs}
(\gsponsor[id=GS4,sponsor-id=100007388]{Compagnia di San
Paolo}/\gsponsor[id=GS5]{Universit\`a di Torino}), \emph{Sviluppo e analisi
di processi\break Markoviani e non  Markoviani con applicazioni}
(\gsponsor[id=GS3]{Universit\`a di Torino}), and by\break
\gsponsor[id=GS7,sponsor-id=100012740]{\xch{INdAM/GNAMPA}{INDAM--GNAMPA}}.
\end{funding}

%
%

\bibliographystyle{abbrvnat}






\end{document}